\documentclass[12pt]{article}
\usepackage{amssymb,amsmath, amsthm, amscd,ifthen}
\usepackage[T2A]{fontenc}
\usepackage[cp1251]{inputenc}
\usepackage[english,russian]{babel}
\usepackage[dvips]{graphicx}
\usepackage[colorlinks,citecolor=red,linkcolor=blue]{hyperref}
\usepackage{indentfirst}
\newtheorem{thm}{Theorem}
\newtheorem{lem}{Lemma}
\newtheorem{cor}{Corollary}

\newtheorem{rem}{Remark}
\newtheorem{prb}{Problem}
\newenvironment{pr1}{\noindent{\it Proof.\,}}{\hfill$\Box$}

\title{Distance graphs with large chromatic number and arbitrary girth}
\author{Andrey Kupavskii}

\begin{document}
\sloppy
\date{}
\maketitle
\renewcommand{\abstractname}{Abstract}

\begin{abstract} In this article we consider a problem related to two famous combinatorial topics. One of them concerns
the chromatic number of the space. The other deals with graphs having big girth (the length of the shortest cycle) and large chromatic number. Namely, we  prove that for any $l\in \mathbb{N}$ there exists a sequence of distance graphs in $\mathbb{R}^n$ with girth at least $l$ and the chromatic number equal to  $(c+\bar{o}(1))^n$ with $c>1$.
\end{abstract}
\section{Introduction}
\subsection {History and related problems}
In this article we study distance graphs  (see \cite{BMP}) of a certain type. Fix some  $a>0$. We say that $G=(V,E)$ is an \textit{$a$-distance graph} in $\mathbb{R}^n$, if $V$ is a subset of $\mathbb{R}^n$ \hypertarget{defdist} and

 $$
E\subseteq\{ \{{\bf x},{\bf y}\}:\ {\bf x},{\bf y} \in V, \,\, \, |{\bf x}-{\bf y}|=a\}.
$$

\begin{rem}If we consider an $a$-distance graph in $\mathbb{R}^n$, then we can apply homothety and transform it into a $1$-distance graph (which is also called unit distance graph). So we won't distinguish $a$-distance graphs for different $a$.
\end{rem}

Such graphs arise naturally in the context of the problem of finding the chromatic number of the space. This famous question was posed by Nelson in 1950: what is the minimum number $\chi(\mathbb{R}^2)$ of colors  needed to color all points of the plane so that no two points at distance one receive the same color? Although this question doesn't sound too difficult, it hasn't got an answer yet. The best we know is that $4\le \chi(\mathbb{R}^2)  \le 7.$
One may ask the same question for higher-dimensional spaces. Here is the formal definition (see \cite{BP}):

$$\chi(\mathbb{R}^n)=\min\{m\in\mathbb{N}: \mathbb{R}^n = H_1\cup \ldots\cup H_m: \\
  \forall i, \forall \mathbf{x},\mathbf{y} \in H_i\ \ |\mathbf{x}-\mathbf{y}|\neq 1\}.$$

  There are quite a few results about this quantity (see the surveys \cite{Rai1}, \cite{Szek} and also \cite{Kup2}), e.g. there are nontrivial lower bounds for the value of
  $\chi(\mathbb{R}^n)$, $n\le 24$. We will be interested in the behavior of $\chi(\mathbb{R}^n)$ as $n\to \infty.$
The following asymptotic lower and upper bounds are due to  A. Raigorodskii \cite{Rai1} and D. Larman, C. Rogers \cite{LR} respectively:

\begin{center}$(\zeta_{low} + o(1))^n\le \chi(\mathbb{R}^n)\le (3+o(1))^n $,\  where $\zeta_{low} = 1.239\dots$
\end{center}

The connection between the chromatic number of the space $\mathbb{R}^n$ and distance graphs in $\mathbb{R}^n$ is intimate. On the one hand,  it follows from the definitions that for any such distance graph $G$, $\chi(G)\le \chi(\mathbb{R}^n).$ On the other hand,  N.G. de Bruijn and P. Erd\H{o}s \cite{ED} proved that there exists a distance graph $G'$ in $\mathbb{R}^n$ with finite number of vertices  such that $\chi(G') = \chi(\mathbb{R}^n).$

The second question that lies at the basis of this article is the following. Can we construct graphs with arbitrarily large chromatic number and arbitrary girth (the length of the shortest cycle)? The positive answer to this question was given by P. Erd\H{o}s \cite{Erd}. He proved that such graphs exist, although the proof was probabilistic, so there was no explicit construction. Later, L. Lov\'asz \cite{Lov} managed to construct such graphs.

It is natural to ask how big can the chromatic number of a distance graph be if we additionally require that the graph has no cliques (complete subgraphs) or cycles of fixed size. The question, whether there is a distance graph in the plane with chromatic number 4 and without triangles (which are both cliques of size 3 and cycles of length 3) was asked by P. Erd\H{o}s \cite{Erd2}. It was answered positively. Moreover, P. O'Donnell (\cite{Don}, \cite{Don2}) proved that for any $l\in \mathbb{N}$ there exists a distance graph in the plane with chromatic number 4 and  girth greater than $l$.

We consider the following three families of distance graphs in $\mathbb{R}^n$: $\mathcal{C}(n,k)$ is the family of all distance graphs that do not contain complete subgraphs of size $k$; $\mathcal{G}_{odd}(n,k)$ is the family of all distance graphs that do not contain odd cycles of length $\le k$; $\mathcal{G}(n,k)$ is the family of all distance graphs that do not contain cycles of length $\le k$. We obviously have the following inclusion: $\mathcal{G}(n,l)\subset\mathcal{G}_{odd}(n,l)\subset \mathcal{C}(n,k)\subset \mathcal{C}(n,k'),$ where $l\ge 3$ and $  k'\ge k\ge 3.$

\hypertarget{defch} Now we define the following quantities:

$$\zeta_k = \liminf_{n\to\infty}\max_{\underline{G}\in \mathcal{C}(n,k)}(\chi(G))^{1/n},$$
$$\xi^{odd}_k = \liminf_{n\to\infty}\max_{\underline{G}\in \mathcal{G}_{odd}(n,k)}(\chi(G))^{1/n},$$
$$\xi_k = \liminf_{n\to\infty}\max_{\underline{G}\in \mathcal{G}(n,k)}(\chi(G))^{1/n}.$$

For example, the bound $\zeta_k\ge 1.1$ means that there exists a sequence of distance graphs $G_{n}\subset \mathbb{R}^n$, such that $\chi(G_n)\ge (1.1+\bar{o}(1))^n$ and none of $G_i$ contains cliques of size $k$.

The values $\zeta_k$ and $\xi^{odd}_k$ were considered in several papers (see \cite{Rai3}). The most accurate estimates on $\zeta_k$ are due to A. Kupavskii \cite{Kup3} (see also \cite{Rar}, where both  $\xi^{odd}_k$ and $\zeta_k$ were considered).

There are two approaches to estimate the quantity $\zeta_k$. The first one is probabilistic, so we don't obtain an explicit graph. However, using this technique one can see that $\zeta_k\ge c_k,$ where $c_k>1$ and $\lim_{k\to \infty}c_k=\zeta_{low}$. In \cite{Rar} this method gave nontrivial bounds only for $k\ge 5.$ A refinement of this method suggested in \cite{Kup3} works for $k\ge 3.$

The second approach is in some sense code-theoretic. It provides us with explicit constructions of such graphs and it works for $k\ge 3$. Moreover, it gives much better bounds for small $k$. But as $k$ grows, this method becomes worse than the probabilistic one, and the bounds tend to some constant that is  significantly smaller than $\zeta_{low}$.

A way to obtain bounds on $\xi^{odd}_k, k\ge 5,$ is also code-theoretic.  In \cite{Rar} it was proved that for any fixed $k$ we have $\xi^{odd}_k>1.$

Bounds on both $\zeta_k$ and $\xi^{odd}_k$ that are slightly weaker than the code-theoretic ones can be derived from simple geometric observations (see \cite{Kup3}).

\subsection{Main Result}
In the previous subsection we considered values $\zeta_k, \xi_k^{odd},
\xi_k.$ For the first two values we mentioned some non-trivial bounds.
However,  both  previous probabilistic and code-theoretic approaches failed to provide any estimate for $\xi_k$, $k\ge 4$. The main result of this article is the following

\begin{thm}\label{th} For any fixed $k\ge 3$ we have $\xi_k\ge 1+\delta,$ where $\delta=\delta(k)$ is a positive constant that depends only on $k$.
\end{thm}

We say that a graph $H$ is a \textit{forest} if it doesn't contain cycles. Consider a finite family $\mathcal{H}  = \{H_1,\ldots, H_m\}$ of graphs. Let $\mathcal{G}(n,\mathcal{H})$ be the family of all distance graphs in $\mathbb{R}^n$ that do not contain any of $H_i\in \mathcal{H}$ as a subgraph. We define the quantity $\xi(\mathcal{H})$ as above:
$$\xi(\mathcal{H}) = \liminf_{n\to\infty}\max_{\underline{G}\in \mathcal{G}(n,\mathcal{H})}(\chi(G))^{1/n}.$$

Theorem \ref{th} provides us with the following appealing corollary:

\begin{cor} For any finite family $\mathcal{H}$ of graphs such that no $H_i\in \mathcal {H}$ is a forest we have
$\xi(\mathcal{H}) \ge 1+\delta,$ where $\delta=\delta(\mathcal{H})$ is a positive constant that depends only on $\mathcal{H}$.
\end{cor}
\begin{pr1}
The proof is immediate. Let $l_i$ be the length of the shortest cycle in $H_i$, where $H_i\in \mathcal{H}$. Then $\xi(\mathcal{H}) \ge \max_{i} \xi_{l_i}\ge 1+\delta,$ where $\delta$ depends only on $\mathcal{H}$.
\end{pr1}

Unfortunately, Theorem \ref{th} says nothing about the family of graphs, on which this bound can be attained. So it is natural to raise the following problem:
\begin{prb} Prove Theorem \ref{th} using an explicit construction.
\end{prb}
The second disadvantage of the method used in this article is the following. The graphs that finally can be obtained are not necessarily \textit{complete distance graphs}, i.e. in the definition of their set of edges we have the strict inclusion (not all possible edges are drawn). Here is another question:

\begin{prb} Prove Theorem \ref{th} using complete distance graphs.
\end{prb}

The rest of the article is organized as follows. In Section \ref{prel} we will give necessary definitions and state auxiliary results. In Section \ref{proof} we will give the proof of Theorem \ref{th}.

\section{Proof of Theorem \ref{th}}\label{sec2}
\subsection{Preliminaries}\label{prel}

As a basis of our construction we will take a family $\mathcal{G}=\{G_{4i}: i\in \mathbb{N} \}$ of distance graphs, where $G_{4n}=(V_{4n},E_{4n}),$ and
$$
V_{4n} = \{{\bf x} = (x_1, \dots, x_{4n}): ~ x_i \in \{0,1\}, ~ x_1 + \ldots + x_{4n} = 2n\}, ~~~
$$
$$
E_{4n} = \{\{{\bf x},{\bf y}\}: ~ ({\bf x},{\bf y}) = n\}.
$$
Here $(,)$ denotes the Euclidean scalar product.
In the next subsection we will prove that for any $k\in \mathbb{N}$ there exists a family of graphs $H_{4i}$ such that for each $i$ the graph $H_{4i}$ is a subgraph of $G_{4i}$ and $H_{4i}$ has girth greater than $k$. Moreover, $\chi(H_{4i})=(c+\delta(4i))^{4i},$ where $c>1$ and $\delta(i)\to 0$ as $i\to \infty.$ This is all we need to prove since in any dimension of the form $4i+j, j=1,2,3$ we can consider a plane of codimension $j$ and embed an isometric copy of $H_{4i}$ there. As a result we obtain a sequence of graphs with desired properties in all dimensions.

It is easy to see that $|V_{4n}| = {4n\choose 2n} = (2+\bar{o}(1))^{4n}$ and $|E_{4n}|={4n\choose 2n}{2n\choose n}^2 = (4+\bar{o}(1))^{4n}.$ We will use the following result from the paper \cite{FR}:

\begin{thm}\label{F} For any $\epsilon>0$ there exists $\delta>0$ such that for any subset $S$ of $V_{4n},$ $|S|\ge (2-\delta)^{4n},$ the number of edges in $S$ (the cardinality of $E_{4n}|_S$) is greater than $(4-\epsilon)^{4n}.$
\end{thm}

\begin{rem} We do not give any numerical bounds for $\xi_k$ since they are very difficult to derive. The reason is that we use Theorem \ref{F}, in which there is no explicit dependency between $\epsilon$ and $\delta.$
\end{rem}

We will also need Lov\'asz Local Lemma (see \cite{AS}):

\begin{thm} \label{lll1} Let $A_1,\ldots, A_m$ be events in an arbitrary probability space and $J(1),\ldots,J(m)$ be subsets of $\{1,\ldots,m\}.$ Suppose there are real numbers $\gamma_i$ such that $ 0<\gamma_i<1,\ i=1,\ldots, m.$ Suppose the following conditions hold:
\begin{enumerate}
\item $A_i$ is independent of algebra generated by $\{A_j, j\not \in J(i)\cup \{i\}\}.$
\item $\mathrm{P}(A_i)\le\gamma_i\prod_{j\in J(i)}(1-\gamma_j).$
\end{enumerate}
Then \ \ \ \ $\mathrm{P}\left(\bigwedge_{i=1}^{m}\overline{A_i}\right)\ge \prod_{i=1}^m(1-\gamma_i)>0.$
\end{thm}

We will use the following version of local lemma (see \cite{Bol}):

\begin{lem}\label{lll} Let $A_1,\ldots, A_m$ and $J(1),\ldots,J(m)$ be as in Theorem \ref{lll1}. Suppose there are real numbers $\delta_i$ such that $ 0<\delta_i\mathrm{P}(A_i)<0.69,\ i=1,\ldots, m.$ Suppose the following condition holds:
\begin{equation}\label{con}
\ln \delta_i\ge\sum_{j\in J(i)}2\delta_j\mathrm{P}(A_j).
\end{equation}
Then \ \ \ \ $\mathrm{P}\left(\bigwedge_{i=1}^{m}\overline{A_i}\right)\ge \prod_{i=1}^m(1-\delta_i\mathrm{P}(A_i))>0.$
\end{lem}
\begin{pr1} This form of local lemma is easy to derive from Theorem \ref{lll1}. We just need to verify that  the inequality 2 from Theorem \ref{lll1} follows from the inequality (\ref{con}). Indeed, we have the following inequality: $$\ln \delta_i\ge\sum_{j\in J(i)}2\delta_j\mathrm{P}(A_j)\ge \sum_{j\in J(i)}-\ln (1-\delta_j\mathrm{P}(A_j)),$$ since $\ln (1-t)\ge -t-t^2\ge -2t$ for $0<t<0.69$ (see \cite{Bol}). We take an exponent of both sides:
$$\delta_i\ge \prod_{j\in J(i)}(1-\delta_j\mathrm{P}(A_j))^{-1}.$$

Finally, we substitute $\delta_i=\gamma_i/\mathrm{P}(A_i)$. 
\end{pr1}

Recall that \textit{the independence number} $\alpha(G)$ of a graph $G =(V,E)$ is the size of a maximum set $S\subset V$ such that for any $v,w\in S$ we have $\{v,w\}\notin E$.
\subsection{Proof of Theorem \ref{th}}\label{proof}

Fix natural numbers $k\ge 3$ and $n$. Let $\gamma\in (0,1)$ be a constant that will be defined later, and set $p=\gamma^{4n}$.
Consider a random subgraph $G$ of the graph $G_{4n}$ in which all edges are chosen independently and uniformly with the probability of each edge to occur equal to $p$. Namely, we have the probability space  $ (\Omega_{4n},{\cal B}_{4n},P_{4n}) $, where
$$
\Omega_{4n} =
\{G=( V_{4n},E), ~ E \subseteq E_{4n}\}, ~~
{\cal B}_{4n} = 2^{\Omega_{4n}}, ~~ P_{4n}(G) = p^{|E|}(1-p)^{|E_{4n}|-|E|} ~~~
{\rm for } ~~ G = (V_{4n},E).
$$

Denote $N = |V_{4n}|.$ We define two families of events on $ \Omega_{4n} $. Firstly, for some $l$ we enumerate all $l$-element subsets of $V_{4n}$ and introduce the events

$$X_i=\{i\text{th $l$-element subset is independent}\},\ \ i=1,\ldots,C_N^l.$$

Secondly, for each $s=3,\ldots, k$ we enumerate all (labeled) cycles of length $s$ in $G_{4n}$ and introduce the events

$$Y_j^s=\{j\text{th $s$-tuple is an $s$-cycle}\},\ \ j=1,\ldots, c_s(G_{4n}),$$

where $c_s(G_{4n})$ is the number of labeled $s$-cycles in $G_{4n}$.

Take $ l = (2-\delta)^{4n} $, where $\delta>0$ is again some constant that will be defined later. The statement of Theorem \ref{th} will follow from the inequality

\begin{equation}\label{prob1}
\mathrm{P}\left(\bigwedge_{i=1}^{C_N^l}\overline{X_i}\wedge\bigwedge_{s=3}^k
\left(\bigwedge_{j=1}^{c_s(G_{4n})}\overline{Y^s_j}\right)\right)>0.
\end{equation}

Indeed, we obtain from (\ref{prob1}) that there exists a subgraph $G'$ in $G_{4n}$ such that it  does not contain cycles of length $\le k$ and at the same time  $\alpha(G')\le l.$ The above means that $G'\in \mathcal{G}(4n,k)$  and
$$\chi(G') \ge \frac{N}{l}=\left(\frac 2 {2-\delta}+\bar{o}(1)\right)^{4n}=\left(1+\delta'+\bar{o}(1)\right)^{4n},$$
where $\delta'$ is a positive constant which will be seen to depend only on $k$.

To prove (\ref{prob1}) we shall use Lemma \ref{lll}. But before we apply it we have to estimate the probabilities of the events $X_i,$ $Y^s_i$.

We start with $ X_i $. Put $a_i=|E(G_{4n}|_{W_i})|,$ where $W_i$ is an $i$-th $l$-element subset of $V_{4n}$. In other words, $a_i$ is the number of edges in $W_i$ in the graph $G_{4n}$.
Then $$\mathrm{P}(X_i)=(1-p)^{a_i}\le e^{-pa_i}=e^{-\gamma^{4n}a_i}.$$

Fix $\epsilon=\epsilon(k)>0$ and choose $\delta=\delta(\epsilon,k)$ as in Theorem \ref{F}. Then obviously $a_i\ge (4-\epsilon)^{4n}.$

We go on to $Y^s_i$. It is easy to see that for each $s$-tuple $Q^s_i$ we have $\mathrm{P}(Y^s_i)=p^{s}=\gamma^{4ns}.$

We also need to analyze the dependencies between the events.

For each $X_i$ let us estimate the number of $Y^s_j$ on which it may depend. Note that if $X_i$ and $Y^s_j$ are dependent then the corresponding sets $W_i$ and $Q^s_j$ must have a common edge. The number of ways to make an $s$-cycle out of a fixed edge is not bigger than $2^{(s-2)4n}$. Thus the number of $Y_j^s$ on which  $X_i$ depends does not exceed  $a_i 2^{(s-2)4n}$.

For all $i,s,$ each $Y^s_i$ and $X_i$ depend on not more than $C_N^l$ events  $X_j.$

It remains to estimate for each $Y^{s_1}_i$ the number of the events $Y^{s_2}_j$ on which it depends. If they are dependent, $Q^{s_1}_i$ and $Q^{s_2}_j$ must have a common edge. Then it is easy to see that the number of such events does not exceed $s_1 2^{4n(s_2-2)}=2^{4n(s_2-2)(1+\bar{o}(1))}.$

For each event $E\in\{X_i, Y^s_i\}$ we split the set $J(E)$ (see Lemma \ref{lll}) into parts. First one ($J^x(E)$) contains all events of the type $X_j.$ The other parts ($J^y_s(E)$) consist of the events of the type $Y^s_j$.
We want to apply Lemma \ref{lll}, so we rewrite the conditions (\ref{con}) for our events:
\begin{equation}\label{sys1}\left\{
\begin{array}{l}
(X_i)\ \ \ \ \ \ln\delta^x_i\ge 2 \sum_{j\in J^x(X_i)}\delta^x_je^{-\gamma^{4n}a_j} + 2\sum_{s=3}^k \sum_{j\in J^y_s(X_i)}\delta^y_j(s) \gamma^{4ns},\\
(Y^{s_1}_i)\ \ \ \ \ \ln \delta^y_i(s_1)\ge 2 \sum_{j\in J^x(Y^{s_1}_i)}\delta^x_je^{-\gamma^{4n}a_j} + 2 \sum_{s=3}^k\sum_{j\in J^y_s(Y_i^{s_1})}\delta^y_j(s) \gamma^{4ns}.
\end{array} \right.
\end{equation}

Fix a constant $f=f(\epsilon, \delta, k)>0,$ which will be defined later, and put $\delta^y_i(s)=e, \delta^x_i= e^{\gamma^{4n(1+f)}a_i}.$ It is easy to see that for sufficiently large $n$ we have $0<\delta^x_i\mathrm{P}(X_i)<0.69, \ 0<\delta^y_i(s)\mathrm{P}(Y_i^s)<0.69$ (see Lemma \ref{lll}). Then for any $j$

$$\delta^x_je^{-\gamma^{4n}a_j}=e^{\gamma^{4n(1+f)}a_j-\gamma^{4n}a_j}=e^{-(\gamma-\bar{o}(1))^{4n}a_j}\le
e^{-(\gamma-\bar{o}(1))^{4n} (4-\epsilon)^{4n}}= e^{-\bigl((4-\epsilon)\gamma-\bar{o}(1)\bigr)^{4n} }.$$

We also have
$$C_N^l \le \left(\frac{eN}{l}\right)^l \le\left(\frac 2{2-\delta}+\bar{o}(1)\right)^{4n\left(2-\delta\right)^{4n}}= e^{\left(2-\delta+\bar{o}(1)\right)^{4n}}.$$
Thus for any $i$
$$\sum_{j\in J^x(X_i)}\delta^x_je^{-\gamma^{4n}a_j}\le \sum_{j\in J^x(X_i)}e^{-\bigl((4-\epsilon)\gamma-\bar{o}(1)\bigr)^{4n} }\le C_N^l e^{-\bigl((4-\epsilon)\gamma-\bar{o}(1)\bigr)^{4n} }\le e^{\left(2-\delta+\bar{o}(1)\right)^{4n}-\bigl((4-\epsilon)\gamma-\bar{o}(1)\bigr)^{4n} }=\bar{o}(1),$$

if \begin{equation}\label{con1} \gamma > \frac{2-\delta}{4-\epsilon}.\end{equation}

 Similarly, if (\ref{con1}) holds, then
$$\sum_{j\in J^x(Y^s_i)}\delta^x_je^{-\gamma^{4n}a_j}=\bar{o}(1).$$

Thereby, if we suppose that (\ref{con1}) holds, then the inequalities (\ref{sys1}) will follow from the system

\begin{equation}\label{sys2}\left\{
\begin{array}{l}
(X_i)\ \ \ \ \ \gamma^{4n(1+f)}a_i\ge  2\sum_{s=3}^k e a_i 2^{4n(s-2)(1+\bar{o}(1))}\gamma^{4ns},\\
(Y^{s_1}_i)\ \ \ \ \ 1\ge  2\sum_{s=3}^k e 2^{4n(s-2)(1+\bar{o}(1))}\gamma^{4ns}.
\end{array} \right.
\end{equation}

One can see that since $\gamma<1$, both inequalities of (\ref{sys2}) are consequences of the following. For any function $g(n)=\bar{o}(1),$ any $s=3,\ldots, k$ and all sufficiently big $n$ should hold
\begin{equation}\label{eq}
2^{4n(s-2)(1+g(n))}\gamma^{4n(s-1-f)} = \bar{o}(1).
\end{equation}

In turn, to prove this it is enough to check the inequality

\begin{equation}\label{eq2}
s-2+ (s-1-f)\log_2 \gamma<0
\end{equation}

for $s=3,\ldots,k.$ Any $\gamma,$

\begin{equation}\label{con2}
0<\gamma<2^{-\frac{k-2}{k-1-f}},
\end{equation}
satisfies (\ref{eq2}), and also the system (\ref{sys2}). Therefore, the system (\ref{sys1}) is satisfied if both (\ref{con1}) and (\ref{con2}) hold:

\begin{equation}\label{con3}
\frac{2-\delta}{4-\epsilon}<\gamma<2^{-\frac{k-2}{k-1-f}}.
\end{equation}

Lastly, we can choose the parameters. We choose $\epsilon = \epsilon(k)$ so that $\frac{2}{4-\epsilon}<2^{-\frac{k-2}{k-1}}.$ Then we choose $f = f(\epsilon,\delta,k)$ small enough so that $\frac{2-\delta}{4-\epsilon}<2^{-\frac{k-2}{k-1-f}}.$ Finally, we choose $\gamma =\gamma(f,\epsilon,\delta,k)\in \left(\frac{2-\delta}{4-\epsilon},2^{-\frac{k-2}{k-1-f}}\right)$.

We have verified all the conditions of Lemma \ref{lll}, hence the inequality (\ref{prob1}) holds and Theorem \ref{th} is proved.

\renewcommand{\refname}{References}


\begin{thebibliography} {20}
\bibitem{AS} N. Alon, J.H. Spencer, {\it The probabilistic method}, New York: Wiley-Interscience, 2000.

\bibitem{BP}
M. Benda, M. Perles, {\it Colorings of metric spaces}, Geombinatorics, 9 (2000), pp. 113-126.

\bibitem{Bol} B. Bollob\'as, {\it Random Graphs}, Cambridge Univ. Press, Second Edition, 2001.

\bibitem{BMP} P. Brass, W. Moser, J. Pach, {\it Research problems in
discrete geometry}, Springer, 2005.

\bibitem {ED} N. G. de Bruijn, P. Erd\H{o}s, {\it A colour problem for infinite graphs and a problem in the theory of relations}, Nederl. Akad. Wetensch. Proc. Ser. A, 54 (1951),  pp. 371–373.

\bibitem{Rar} E. Demechin, A. Raigorodskii, O. Rubanov, {\it Distance graphs with big chromatic number that do not contain cliques or cycles of fixed size}, Sbornik: Mathematics, to appear.

\bibitem{Don} P. O'Donnell, {\it Arbitrary girth, 4-chromatic unit distance graphs in the plane I. Graph embedding}, Geombinatorics,  9 (2000), pp. 180-193.

\bibitem{Don2} P. O'Donnell, {\it Arbitrary girth, 4-chromatic unit distance graphs in the plane II. Graph description}, Geombinatorics,  9 (2000), pp. 145-152.

\bibitem{Erd} P. Erd\H{o}s, {\it Graph theory and probability}, Canad. J. Math., 11 (1959), pp. 34-38.

\bibitem{Erd2} P. Erd\H{o}s, {\it Unsolved Problems}, Congres Numerantium XV -- Proceedings of the 5th British Comb. Conf. 1975, (1976), p. 681.

\bibitem{FR} P. Frankl, V. R\"odl, {\it Forbidden intersections}, Trans. of Amer. Math. Soc., 300 (1987), N1, pp. 259-286.

\bibitem{Lov} L. Lov\'asz, {\it On Chromatic Number of Finite Set-Systems,} Acta Math. Acad. Sci. Hungar. 19 (1968), pp. 59-67.

\bibitem{Kup2}
A. Kupavskii, {\it On the colouring of spheres embedded in $\mathbb{R}^n$}, Sbornik: Mathematics, 202 (2011), N6,  pp. 859-886.

\bibitem{Kup3} A. Kupavskii, {\it Explicit and probabilistic constructions of distance graphs with small clique numbers and big chromatic numbers}, Izvestiya mathematics, submitted.

\bibitem{LR} D.G. Larman, C.A. Rogers, {\it The realization
of distances within sets in Euclidean space}, Mathematika, 19 (1972), pp. 1-24.


\bibitem{Rai1} A. Raigorodskii, {\it Borsuk's problem and the chromatic numbers of some metric spaces}, Russian Mathematical Surveys,  56(2001), N1, pp. 103–139.

\bibitem{Rai3} A. Raigorodskii, {\it On distance graphs with large chromatic number but without large simplices}, Russian Mathematical Surveys,  62 (2007), N6, pp. 1224–1225

\bibitem{Szek} L.A. Sz\'ekely, {\it Erd\H{o}s on unit distances
and the Szemer\'edi - Trotter theorems}, Paul Erd\H{o}s and his Mathematics, Bolyai Series Budapest, J. Bolyai Math.
Soc., Springer, 11 (2002), pp. 649-666.



\end{thebibliography}
\end{document}